\documentclass[12pt]{amsart}
\usepackage[margin=1in]{geometry}
\usepackage{amsthm,amssymb, mathtools}
\usepackage{amsmath,xypic}
\usepackage{graphicx}
\usepackage[numbers]{natbib}
\usepackage{fancyhdr} % used by RCS
\newcommand{\be}{\begin{equation}}
\newcommand{\ee}{\end{equation}}
\newcommand{\bes}{\begin{equation*}}
\newcommand{\ees}{\end{equation*}}
\newcommand{\bea}{\begin{eqnarray}}
\newcommand{\eea}{\end{eqnarray}}
\newcommand{\beas}{\begin{eqnarray}}
\newcommand{\eeas}{\end{eqnarray}}
\newcommand{\ben}{\begin{note}}
\newcommand{\een}{\end{note}}
\newcommand{\bexl}{\vskip0.1em\noindent\hrulefill\vskip1em\begin{ExerciseList}}
\newcommand{\eexl}{\end{ExerciseList}\hrulefill}

\newcommand{\bthm}{\begin{theorem}}
\newcommand{\ethm}{\end{theorem}}
\newcommand{\bpro}{\begin{prop}}
\newcommand{\epro}{\end{prop}}
\newcommand{\bcor}{\begin{corollary}}
\newcommand{\ecor}{\end{corollary}}
\newcommand{\bcon}{\begin{conjecture}}
\newcommand{\econ}{\end{conjecture}}
\newcommand{\bp}{\begin{proof}}
\newcommand{\ep}{\end{proof}}
\newcommand{\blem}{\begin{lemma}}
\newcommand{\elem}{\end{lemma}}
\newcommand{\bn}{\begin{note}}
\newcommand{\en}{\end{note}}
\newcommand{\benum}{\begin{enumerate}}
\newcommand{\eenum}{\end{enumerate}}
\newcommand{\bed}{\begin{defn}}
\newcommand{\eed}{\end{defn}}
\newcommand{\br}{\begin{remark}}
\newcommand{\er}{\end{remark}}

\newcommand{\lab}[1]{\label{#1}}

\newtheorem{theorem}[equation]{Theorem}      % (If you want theorem numbered
\newtheorem{lemma}[equation]{Lemma}          %
\newtheorem{corollary}[equation]{Corollary}  %       goes for lemmas, etc.)
\newtheorem{proposition}[equation]{Proposition}

\newtheorem{conj}[equation]{Conjecture}
\newtheorem*{conjecture*}{Conjecture}
\newtheorem*{theorem*}{Theorem}

\theoremstyle{definition}
\newtheorem{conjecture}[equation]{Conjecture}

\theoremstyle{definition}
\newtheorem{defn}[equation]{Definition}
\theoremstyle{remark}

\theoremstyle{definition}
\newtheorem{remark}[equation]{Remark}

\numberwithin{equation}{section}
\newcommand{\N}{\mathbb{N}}

\newcommand{\Q}{{\mathbb Q}}

\renewcommand{\int}{\operatorname{int}}

\renewcommand{\P}{{\mathbb P}}

\newcommand{\QQ}{\mathbb{Q}}
\newcommand{\ZZ}{\mathbb{Z}}

\newcommand{\NN}{\mathbb{N}}

\newcommand{\Oc}{\mathcal{O}}

\newcommand{\Cond}{\mathrm{Cond}}
\newcommand{\CondP}{\mathrm{Cond}_{\P}}
\newcommand{\CondAP}{\mathrm{Cond}_{\AP}}
\newcommand{\Condss}{\mathrm{Cond}_{\textrm{ss}}}

\renewcommand{\P}{\mathcal{P}} % for Primes
\newcommand{\AP}{\mathcal{AP}} % for almost Primes
\renewcommand{\bpro}{\begin{proposition}}
\renewcommand{\epro}{\end{proposition}}
\begin{document}

\title[]{Asymptotics of conductors of elliptic curves over $\mathbb{Q}$}%
\author{Sean Howe}
\address{Dept. of Mathematics, 
University of Chicago, 5734 S. University Avenue, Chicago IL 60637}
\email{seanpkh@gmail.com}

\author{Kirti Joshi}%
\address{Math. department, University of Arizona, 617 N Santa Rita, Tucson AZ
85721-0089, USA.} \email{kirti@math.arizona.edu}
%\date{Preliminary Version: \today}

\thanks{}%
\subjclass{}%
\keywords{}%

%\date{}%
%\dedicatory{}%
%\commby{}%

% ----------------------------------------------------------------
\begin{abstract}
In this note we study numbers which occur as conductors of elliptic curves over $\Q$. We show, by constructing families of elliptic curves with quadratic discriminant and invoking a theorem of Iwaniec, that this set contains infinitely many almost primes. We show, assuming a strong version of the Cohen-Lenstra heuristics, that the set of prime conductors has an explicitly bounded density in the set of primes. Studying the Cremona and Stein-Watkins databases of elliptic curves we conjecture that the set of conductors should be of density zero in the set of natural numbers and that the set of prime conductors should be of density zero in the set of prime numbers.
\end{abstract}
\maketitle

%-----------------------------------------------------------------

% ----------------------------------------------------------------

\let\QQ=\Q
\let\NN=\N

\tableofcontents

\section{Introduction}\label{intro}
Let $n\geq 1$ be a natural number. We say that \textit{$n$ is a conductor} or more precisely \textit{an elliptic conductor} if there exists an elliptic curve $E/\Q$ with conductor equal to $n$. If there is no elliptic curve over $\Q$ of conductor $n$, then we say that $n$ is a \textit{non-conductor}. 
By the Modularity Theorem, the existence of an elliptic curve of conductor $n$ is equivalent to the existence of a weight 2 eigenform with rational Hecke eigenvalues in spaces of newforms of level $\Gamma_0(n)$. These spaces are effectively computable, and using this, or by other methods, one knows, e.g., that the natural numbers $1\leq n\leq 10$, are not conductors. 

For any set of natural numbers $S$, let $\Cond_S$ be the set of all elliptic conductors in $S$. In this work, we study asymptotic properties of $\Cond_S$ for the following sets $S\subset\N$: $S=\N$, $S=\P$ the set of primes, $S=\AP$ the set of almost primes ($n$ is an almost prime if $n$ has at most 2 prime factors, counted with multiplicity) and $S=ss$ the set of square-free integers. Note $n\in\Cond_{ss}$ if and only if $n$ is a conductor and $n$ is square-free, or equivalently there is a semi-stable elliptic curve of conductor $n$, thus the notation $S=ss$. Clearly 
 \be
 \Cond\supset\CondAP\supset\CondP
 \ee
 and 
 \be 
    \Cond\supset\Condss\supset\CondP.
 \ee  
For $X\geq 1$, and for $T \subset \N$ we denote $T(X)=T\cap[1,X]$. In each of the above cases, we are interested in the following basic questions:
\begin{enumerate}
\item Is $\Cond_S$ infinite? What about its complement (in $S$)?
\item If $\Cond_S$ is infinite, does $\Cond_S$ have a density (in $S$)?
\item If the density of $\Cond_S$ is zero, what is the asymptotic growth rate of $|\Cond_S(X)|$?
\end{enumerate}
As noted above, $2 \not\in \Cond$, thus $\Cond$ is a proper subset of $\NN$ and $\CondP$  (resp. $\CondAP$) is a proper subset of the set of
primes $\P$ (resp. set of almost primes $\AP$).

Standard results on conductors show that both $\Cond$ and its complement in $\NN$ are infinite (cf. Section \ref{SecCond}). Based on numeric evidence from the Stein-Watkins database, we make the following conjecture:

\begin{conjecture*}[cf. Conjectures \ref{con:zero-density} and \ref{con:prime-zero-density}]
$\Cond$ has density zero in $\NN$, and $\CondP$ has density zero in $\P$. 
\end{conjecture*}

This and related conjectures are discussed further in Section \ref{sec:Conjectures}. 

In Section~\ref{SecCondP}, the bulk of this note, we study prime and almost prime conductors. Our main results are the following:

\begin{theorem*}[cf. Theorems \ref{prop-InfPrimeCondWithTwoTorsion} and \ref{prop-InfPrimeCondNo2Torsion}] \hfill
\begin{enumerate}
\item Assuming the conjecture of Hardy and Littlewood \cite[Conjecture F]{HardyLittlewood-QuadraticConjecture}, the set $\CondP$ is infinite.
\item The sets $\CondAP$ and $\Condss$ are infinite.
\end{enumerate}
\end{theorem*}

\begin{theorem*}[cf.  Theorem \ref{ThmModDensity} and Corollary \ref{CorTotalDensity}] Assuming a Cohen-Lenstra type conjecture (Conjecture \ref{CLPrecise3Div}), the upper density of $\CondP$ in $\P$ is bounded by .69, and in particular, there are infinitely many primes that do not appear as conductors. 
\end{theorem*}

Although there is a considerable classical body of work on conductors concerned with the explicit computation of the elements of sets such as $\CondP$ under a fixed size in order to provide evidence for the Modularity Theorem (cf., e.g., \cite{Ogg-AC2PC}, \cite{Ogg-ACSC}, \cite{Setzer-ECPC}, \cite{Hadano-CECRPO2}), our two theorems are, to our knowledge, the first results on the asymptotic behavior of these sets. 

We now outline the proofs of these two theorems. A classic result of Setzer \cite[Theorem 2]{Setzer-ECPC} states that for $p>17$ there exists an elliptic curve of conductor $p$ with a rational point of order 2 if and only if $p=u^2 + 64$ for an integer $u$, and explicitly constructs all such curves given such a $p$. A well-known conjecture of Hardy and Littlewood \cite{HardyLittlewood-QuadraticConjecture} then implies that there exist infinitely many such primes $p$, and thus that $\CondP$ is infinite. Using a similar argument we give infinite families of primes $p$ appearing as conductors of curves with trivial rational 2-torsion, again conditional on the conjecture of Hardy and Littlewood and each contributing zero to the density of $\CondP$.  Without assuming the conjecture of Hardy and Littlewood we are unable to prove that $\CondP$ is infinite, however, a result of Iwaniec implies unconditionally that the families we construct contain infinitely many curves of distinct almost prime conductors and infinitely many curves with distinct semistable conductors. We note that in order to produce infinitely many semistable conductors divisible by at most two primes the Setzer families do not suffice and it is necessary to use the additional families we construct (cf. Remark \ref{rem:TwoTorsionCond}-(1)). 

Bounding the density of $\CondP$ in $\P$ is more delicate. By the result of Setzer quoted above and a trivial estimate, we are reduced to considering curves with no rational 2-torsion. Then, a classical approach dating back to work of Ogg \cite{Ogg-ACSC} and Hadano \cite{Ogg-ACSC} translates the existence of such a curve into 3-divisibility of class numbers of associated quadratic fields. Thus, to restrict the existence of such curves, it suffices to restrict the existence of quadratic fields with certain class numbers, and we enter the realm of Cohen-Lenstra heuristics.

To apply Cohen-Lenstra in this situation, we formulate a version that applies to quadratic fields $Q(\sqrt{\pm p})$ for $p$ prime and in a specific congruence class mod $8$. The differences from the classical situation \cite{CohenLenstra} are the restriction to primes, the splitting into congruence classes mod 8, and an independence assumption for real and imaginary quadratic statistics. We discuss reasons to believe this conjecture in Remark \ref{rem:CLWhy}. 

The restriction to different congruence classes mod $8$ is a key ingredient in our proof of the upper bound, and indeed the cumulative upper bound we give is built out of distinct upper bounds for each congruence class. Furthermore, although we conjecture that the density of $\Cond_P$ in $\P$ is zero, the precise asymptotics also appear to depend on the congruence class mod $8$ in a manner consistent with the difference in our upper bounds. This is one of the more interesting features we observe in the asymptotic behavior of $\CondP$; cf. also Remark \ref{rem:Families}-(3).

There are natural generalizations of the questions we consider here, which can be phrased either in the language of automorphic forms or Galois representations. For example, we might consider modular forms of higher weights, or eigenforms with coefficients in a larger field, and ask in which levels they appear. However, because our approach depends heavily on the concrete interpretation provided by elliptic curves over $\QQ$, and because most of the numeric data available is on elliptic curves, we do not pursue these generalizations here.

\section{Some facts about conductors}\lab{SecCond}
We begin by summarizing some standard facts about the conductors of elliptic curves (see \cite{Silverman-AOEC}):

\begin{proposition}[see \cite{Silverman-AdvTop}]\label{PropCondProperties}
Let $E/\QQ$ be an elliptic curve and let $N_E = \prod_p p^{f(p)}$ be the conductor of $E$. Then,
\begin{enumerate}
\item $f(p)=0$ if and only if $E$ has good reduction at $p$
\item $f(p)=1$ if and only if $E$ has multiplicative reduction at $p$
\item For $p>3$, $f(p)=2$ if and only if $E$ has additive reduction at $p$
\item For $p=3$, $2\leq f(p) \leq 5$ if and only if $E$ has additive reduction at $p$
\item For $p=2$, $2\leq f(p) \leq 8$ if and only if $E$ has additive reduction at $p$.
\end{enumerate}
\end{proposition}

In particular, $\Cond$ is contained in the set of ninth-power free integers. The set of ninth-power free integers is of density $\frac{1}{\zeta(9)}=0.9979956\cdots$ and thus its complement is infinite (in $\N$) and as no elliptic curve can have conductor, say $p^3$, for a prime $p>3$, it is immediate that the complement of $\Cond$ in the set of ninth-power free integers is also infinite. The set $\Cond$ has further structure coming from quadratic twists. For example:

\begin{proposition}[see  \cite{Silverman-AdvTop}]
If $n \in \Cond$ and $m$ is a squarefree integer with $(m,n)=1$ then  $m^2 n \in \Cond$.
\end{proposition}
\begin{proof}
Let $m'$ be the part of $m$ not divisible by $2$, and take a quadratic twist by $\pm m'$ of a curve of conductor $n$ where the sign chosen depends on whether or not $2 | m$. 
\end{proof}

This shows that $\Cond$ is infinite -- for example, since there is a curve of conductor 11, there is a curve of conductor $11m^2$ for every square-free $m$ coprime to $11$. Alternatively, one can deduce that $\Cond$ is infinite from the result for almost prime conductors (Theorems \ref{prop-InfPrimeCondWithTwoTorsion} and \ref{prop-InfPrimeCondNo2Torsion}), which gives explicit families of curves containing infinitely many curves of almost prime conductor (which, in particular, will not be twists of one another for distinct conductors). 

\begin{remark}
The previous proposition suggests the following definition: a conductor $n$ in $\Cond$ is called primitive if $n$ cannot be written as $m^2 n'$ for a square-free $m > 1$ and $n' \in \Cond$ such that $(m,n')=1$. The set $\Cond$ is then a kind of convolution of the set of primitive conductors and the set of square-free integers, allowing us to reduce some questions about $\Cond$ to questions about primitive conductors. Primitive conductors are related to other sets of conductors that we consider: it is elementary that $\Cond_{ss}$ is contained in the set of primitive conductors and $1 \not \in \Cond$ implies $\Cond_\AP$ is contained in the set of primitive conductors. Since neither $\Cond_\AP$ nor $\Cond_{ss}$ contains the other (e.g., $49 \in \Cond_\AP$ and $30 \in \Cond_{ss}$), both these inclusions are proper. Although primitive conductors are potentially a useful tool in studying $\Cond$, since the majority of this article is focused on prime and almost prime conductors, we do not develop the idea further.  
\end{remark}

Our main tool for analyzing the conductor will be the minimal discriminant. We record here some relations between the discriminant and conductor:

\begin{proposition}[cf. \cite{Silverman-AOEC}, Proposition 5.1] \label{PropCondDisc}
Let $E/\QQ$ be an elliptic curve of conductor $N$, let $\Delta$ be the minimal discriminant of $E$ and $c_4$ the invariant of that name associated to a minimal Weierstrass equation. Let $p$ be a prime. Then
\begin{enumerate}
\item The following are equivalent:
\begin{itemize}
	\item E has good reduction at $p$
	\item $p \not | \Delta$
	\item $p \not | N$
\end{itemize}
\item The following are equivalent:
\begin{itemize}
	\item E has multiplicative reduction at $p$
	\item $p | \Delta$ and $p \not | c_4$
	\item $p$ divides $N$ exactly once
\end{itemize}
\item If $p \geq 5$ and $E$ has additive reduction (equivalently $p^2 | N$), then $p^2 | \Delta$.
\end{enumerate}
\end{proposition}
\begin{proof}
The last item (3) follows from the identity $1728\Delta = c_4 ^3 - c_6 ^2$ and Proposition \ref{PropCondProperties}. 
\end{proof}

\section{Prime and almost prime conductors}\label{SecCondP}
In this section, the bulk of our work, we examine the sets $\CondP$ of prime conductors and $\CondAP$ of almost prime conductors. In Section \ref{subsec-InfPrimeCond} we show that $\CondAP$ is infinite and, conditional on a conjecture of Hardy Littlewood, that $\CondP$ is infinite. To do so we examine explicit families with a rational 2-torsion point (see  Theorem \ref{prop-InfPrimeCondWithTwoTorsion}) and without a rational 2-torsion point (see  Theorem \ref{prop-InfPrimeCondNo2Torsion}) with quadratic discriminant and show that Iwaniec's theorem, resp. the conjecture of Hardy and Littlewood, implies that infinitely many of the curves in both these families have almost prime, resp. prime, conductors.

In Section \ref{subsec:BoundingDensities} we show that the existence of a curve with prime conductor and no trivial rational 2-torsion gives a non-trivial class number condition on the fields $\QQ(\sqrt{\pm p)}$. We then formulate a version of the Cohen-Lenstra heuristics (Conjecture \ref{CLPrecise3Div}) that applies to this situation, and use it to deduce a conditional upper bound on the density of $\Cond_P$ in $P$ (Theorem \ref{ThmModDensity} and Corollary \ref{CorTotalDensity}).

\subsection{Constructing curves with prime and almost prime conductor} \label{subsec-InfPrimeCond}
Setzer \cite[Theorem 2]{Setzer-ECPC} showed that for $p\neq2,\,3,\,17$
there is an elliptic curve of conductor $p$ with a rational 2-torsion
point if and only if $p=u^{2}+64$ for some integer value of $u$.
As part of this construction he considers, for any odd value of $u$,
the elliptic curve given by the Weierstrass equation (where the sign
of $u$ is chosen so that $u\equiv1\mod4$)
\begin{equation}\label{eq-SetzerCurve}
y^{2}+xy=x^{3}+\frac{1}{4}(u-1)x^{2}-x
\end{equation}
which has discriminant $\Delta=u^{2}+64$. When $p=u^2+64$ for $p$ prime, this gives a minimal
model. We will see that in this case the curve has conductor $p$.

A well-known conjecture of Hardy and Littlewood \cite[Conjecture F]{HardyLittlewood-QuadraticConjecture}  (see \cite{BaierAndZhao-PrimesByQuadratics} for a more recent
survey) on primes represented by quadratic polynomials implies that there are infinitely many values of $u$ such that $u^2+64$ is prime.  Concretely, the
conjecture implies that there exists a positive constant $C$ such that
\be
p(x)\sim C\cdot\frac{\sqrt{x}}{\log x}
\ee
where $p(x)$ is the number of primes $q\leq x$ such that $u^2+64=q$
for some integer $u$, and gives a formula for the constant $C$. It is a theorem due to Iwaniec \cite[Section 1]{IwaniecAlmostPrimes} that $u^2 + 64$ represents infinitely many almost primes. In particular, we obtain the following result

\begin{theorem} \label{prop-InfPrimeCondWithTwoTorsion}
\begin{enumerate}
\item Assuming the conjecture of Hardy and Littlewood \cite[Conjecture F]{HardyLittlewood-QuadraticConjecture}, there are infinitely many primes $p$ such that there exists an elliptic curve of conductor $p$ with a rational point of order 2.
\item There are infinitely many almost primes $n$ such that there exists an elliptic curve of conductor $n$ with a rational point of order 2.
\item The density of the set of primes $p$ such that there exists an elliptic curve $E/\QQ$ of conductor $p$ with a rational point of order 2 is zero (in the set of all primes).
\end{enumerate}
\end{theorem}
\begin{proof}
The first statement follows from the discussion above and the fact that any prime of additive reduction greater than or equal to $5$ divides the minimal discriminant to a power greater than $1$ (Proposition \ref{PropCondDisc}). 

For the second statement, we observe that the result of Iwaniec \cite[Section 1]{IwaniecAlmostPrimes} implies there are infinitely many $u$ such that the Setzer curve \ref{eq-SetzerCurve} has almost prime discriminant $\Delta$ -- to apply it in this case, make the change of variables $u=x+1$. Such a $u$ must be odd since otherwise $4 | u^2 + 64$ which also must have another prime factor. Thus, $2$ does not divide $\Delta$. Observe $c_4 = u^2 + 48$ so $\Delta - c_4 \equiv 1 \mod 3$, thus if $3 | \Delta$ the curve has multiplicative reduction at 3 and the power of $3$ dividing the conductor is exactly 1. Since $\Delta$ is almost prime, any prime $p\geq 5$ dividing $\Delta$ is either of multiplicative reduction or $\Delta = p^2$. 

The final statement follows from the prime number theorem and the easy estimate that there are at most $\sqrt{x}$ numbers less than $x$ of the form $u^2+64$. 
\end{proof}

\begin{remark} \label{rem:TwoTorsionCond} We make a few comments on Theorem \ref{prop-InfPrimeCondWithTwoTorsion}:
\begin{enumerate}
\item It is a priori possible that the infinitely many almost prime conductors given by (2) and represented by $u^2+64$ are all of the form $p^2$ for primes $p$. In particular, in contrast to Theorem \ref{prop-InfPrimeCondNo2Torsion}, we cannot use (2) to deduce that $\Cond_{ss}$ is infinite. 
\item One might expect (3) to hold for almost primes as well, but without further analysis we cannot make a similar statement about almost primes because not all curves with a rational two torsion point and almost prime conductor appear in the Setzer families. 
\item We thank the anonymous referee who pointed out the change of coordinates used to apply Iwaniec's theorem in (2).
\end{enumerate}
\end{remark}

In fact, for this kind of argument there is nothing special about curves with a rational 2-torsion point except that those of prime conductor all live in the one-parameter family of Setzer curves (\ref{eq-SetzerCurve}) whose discriminant is an irreducible quadratic polynomial in one variable. It is easy to produce more such families:

\begin{proposition}
\label{prop-quadDiscFamily}Let $a,b\in\mathbb{Z}$, and consider the Weierstrass equation
\begin{equation}
y^{2}+y=x^{3}+ax^{2}+bx+n\label{eq:MinimalModel}
\end{equation}
with discriminant
\be
\Delta(n)=(-432)n^{2}+(-64a^{6}+288a^{2}b-216)n+(-16a^{6}+16a^{4}b^{2}-64b^{3}+72a^{2}b-27).
\ee
If at least one of $a$ and $b$ is not divisible by 3 and $(a^4-3b^3)$ is not a square, then $\Delta(n)$ is an irreducible quadratic polynomial (in $\ZZ[n]$).
\end{proposition}
\begin{proof}
Considering $\Delta(n)$ as a polynomial in $n$,
the leading coefficient -432 is divisible only by 2 and 3, but the
constant coefficient is odd. Since our condition on $a,b\mod3$ implies
that either the coefficient of $n$ or the constant coefficient is
not divisible by 3, the gcd of the coefficients is 1.  Furthermore, a calculation shows that if we consider $\Delta(n)$ as a quadratic polynomial in $n$ then its discriminant is equal to
$(a^4-3b^3)^3$
which by hypothesis is not square and thus $\Delta(n)$ is an irreducible
polynomial in $n$.
\end{proof}

\begin{remark}
The condition that ${a}^{4}-3{b}^{3}$ is a square is satisfied, e.g., for $(a,b)$ equal to
 $(0,1)$, $(1,0)$, $(1,1)$. For a fixed generic value of $b$
(resp. $a$) the curve
\be
y^{2}=a^4 - 3b^3
\ee
is smooth of genus 1, thus by Siegel's theorem there will be at
most finitely many integer values of $a$ (resp. $b$) such that this
quantity is a square.
\end{remark}

\begin{theorem}\label{prop-InfPrimeCondNo2Torsion}
\begin{enumerate}
\item 
Assuming the conjecture of Hardy and Littlewood \cite[Conjecture F]{HardyLittlewood-QuadraticConjecture}, there exist infinitely many primes $p$ such that there exists an elliptic curve $E/\QQ$ of conductor $p$ with trivial rational 2-torsion.
\item There exist infinitely many $n\in \Cond$ such that $n$ is prime or $n=pq$ for $p$ and $q$ distinct odd primes. In particular, $\CondAP$ and $\Condss$ are both infinite sets. 
\end{enumerate}
\end{theorem}
\begin{proof}
For the first item, we observe that for any value of $a$ and $b$ in Proposition \ref{prop-quadDiscFamily}, we have $\Delta(n)\equiv 5\mod8$. In particular, if $\Delta(n)=\pm p$ for $p>17$ prime then, arguing as in Theorem \ref{prop-InfPrimeCondWithTwoTorsion}, we obtain a curve $E$ of conductor $p$ that cannot have a rational point of order 2 (because by the result of Setzer we have $p\equiv 1 \mod 8$ if $E$ has a rational point of order 2). Thus, applying Hardy and Littlewood as before to the discriminant of the family of curves \ref{eq:MinimalModel} with, for instance, $a=0,\, b=1$, we obtain the result. 

For the second item, we again consider the families of curves (\ref{eq:MinimalModel}) of Proposition \ref{prop-quadDiscFamily} for a fixed value of $a$ and $b$, e.g. $a=0$, $b=1$. By Iwaniec's theorem \cite[Section 1]{IwaniecAlmostPrimes}, the discriminant $\Delta(n)$ of these curve takes on infinitely many almost prime values. Because \[\Delta(n) \equiv 5 \mbox{ mod } 8,\] none of these almost-primes is equal to $\pm p^2$ for a prime $p$, or divisible by $2$, and thus the conductor is divisible by two distinct primes $\geq 3$. We observe that for Equation \ref{eq:MinimalModel}, $c_{4}=16(a^{4}-3b)$, which is divisible by 3 if and only if $a$ is. But if $ a\equiv 0\mod 3$,
then $\Delta(n)\equiv 2b^{3} \mod 3$ and since by hypothesis not both of $a$ and $b$ are divisible by 3, we see that $\Delta(n)$ is not divisible by
3. Thus if $3 | \Delta(n)$, $3 \not| c_4$, and the curve has multiplicative reduction at $3$. We conclude by Proposition \ref{PropCondDisc} that the conductor equals $|\Delta(n)|$ and thus is a product of two distinct primes. 
\end{proof}

\begin{remark} \label{rem:Families} We make some comments on Theorem \ref{prop-InfPrimeCondNo2Torsion}:
\begin{enumerate}
\item It may be possible to choose a family so that the curves in (2) have no rational two torsion points, but we have not investigated this further. 
\item Proposition \ref{prop-quadDiscFamily} gives infinitely many families of curves, but in each the the discriminant is quadratic, so we cannot use it to give even a conditional non-trivial lower bound on the density of $\CondP$ in $\P$. 
\item The primes $p \in \CondP$ coming from Theorem \ref{prop-InfPrimeCondWithTwoTorsion}, that is, from the Setzer family (\ref{eq-SetzerCurve}) are all of the form $p \equiv 1 \mod 8$. Similarly, the infinitely many primes with trivial rational two torsion of Theorem \ref{prop-InfPrimeCondWithTwoTorsion} coming from the families of Proposition \ref{prop-quadDiscFamily} are $\equiv 3 \mod 8$ (since $\Delta(n)\equiv 5 \mod 8$ and the leading coefficient is negative). It may be possible to write down elementary families in the same fashion giving prime conductors $\equiv 5,7 \mod 8$, however we have not investigated the question closely, and in fact there is reason to believe it is more difficult. Indeed, in Section \ref{subsec:BoundingDensities} we will see that the existence of an elliptic curve $E$ with minimal discriminant $d=\pm p$ for $p$ prime gives a non-trivial condition on the class number of $\QQ(\sqrt{d})$ except when $E$ has non-trivial rational 2-torsion or $d<0$ and $d\equiv 5\mod 8$, and thus the Setzer curves and the families produced above live precisely in these two exceptional cases. For example, a one-parameter family of curves with positive quadratic discriminant $d \equiv 7 \mod 8$ combined with the conjecture of Hardy and Littlewood would imply the existence of infinitely many primes $p$ such that $Q(\sqrt{\pm p})$ had class number divisible by $3.$

\item By fixing different coefficients and/or considering other families it is easy to obtain one parameter families of elliptic curves with higher degree polynomial discriminants (for example, if we fixed $n$ and $a$ and varied $b$ in the family (\ref{eq:MinimalModel}) of Proposition \ref{prop-quadDiscFamily} then we would obtain a cubic polynomial in $b$ for the discriminant), and one could apply conjectures analogous to those of Hardy and Littlewood to these families. Similarly one could allow all of the coefficients to vary simultaneously and apply conjectures on primes represented by multivariate polynomials. We have focused on the one-variable quadratic case because of the simple examples it provides and because of the results of Iwaniec \cite{IwaniecAlmostPrimes} on almost-primes represented by quadratics.

\item The conjecture of Hardy and Littlewood applied to the polynomial $\Delta(n)$ of Proposition \ref{prop-quadDiscFamily} gives the precise asymptotic estimate
\be
p(x)\sim C\cdot\frac{\sqrt{x}}{\log x}
\ee
where $p(x)$ is the number of primes $q\leq x$ such that $\Delta(n)=q$
for some integer $n$, and gives a formula for the constant $C$,
which in our case will depend on $a$ and $b$. For example, when
$a=b=1$ our theorem applies,
\be
\Delta(n)=-432n^{2}+8n-19,
\ee
and we calculate $C\sim0.063$. When $a=0$ and $b=1$, $C\sim0.162$.

\item 
In \cite{iwaniec98} it was shown that $y^2+x^4=p$ has solutions $(x,y)$ for infinitely many primes $p$. However note that $y^2=p-x^4$ has Complex Multiplication by $\Q(\sqrt{-1})$ and hence it has potentially good reduction at all primes dividing its conductor and it is evident from the equation that the conductor is divisible by $4$ and $p^2$ so this curve does not provide elliptic curves of prime conductor.
\end{enumerate}
\end{remark}

\subsection{Bounding the density of prime conductors}\label{subsec:BoundingDensities}

In this section we bound the density of primes appearing as conductors of elliptic curves, conditional on a Cohen-Lenstra type conjecture (Conjecture \ref{CLPrecise3Div}). Our main results are Theorem \ref{ThmModDensity} and Corollary \ref{CorTotalDensity}.

By Theorem \ref{prop-InfPrimeCondWithTwoTorsion}-(3), to bound the density of prime conductors, it suffices to bound the density of prime conductors appearing as the conductors of elliptic curves with trivial rational 2-torsion. Such a curve for a prime $p$ gives rise to a restriction on the ray class group of conductor (2) of $\QQ(\sqrt{\pm p})$ (Proposition \ref{PropRamTwoTorsion}), which we then analyze using the Cohen-Lenstra heuristics. 

The idea of using the class group of $\QQ(\sqrt{\pm p})$ to study curves of prime conductor is classic and our Proposition \ref{PropRamTwoTorsion} is essentially already contained in work of Ogg \cite{Ogg-AC2PC, Ogg-ACSC} and Setzer \cite{Setzer-ECPC}, and in a more refined form in work of Brumer and Kramer \cite{Brumer-ROEC}. However because it does not take long we reprove it here in order to obtain the precise statement we need. The only novelty in our statement is the replacement of the class group with the ray class group of conductor (2), which allows us to make non-trivial statements in the case of prime conductors $p \not \equiv \pm 1 \mod 8$.

\begin{lemma}\label{LemMultRed}
Let $E/\QQ$ be an elliptic curve, and let $L=\QQ(E[2])$ be the two-torsion field. For $p$ a prime of $\QQ$, denote by $e_p$ and $f_p$ its ramification degree and residue extension degree in $L$. If $p$ is a prime of multiplicative reduction for $E$ then $f_p \cdot e_p | 2$.
\begin{proof}
Tate's $p$-adic uniformization shows that for $p$ a prime of multiplicative reduction, there exists a $\QQ_p$-rational $2$-torsion point, and thus $\QQ_p(E_{\QQ_p}[2]) / \QQ_p$ is an extension of degree $1$ or $2$.
\end{proof}
\end{lemma}

\begin{remark} For $p > 2$ it is not necessary to invoke uniformization -- one can proceed directly by extracting a root of $f$ in a Weierstrass equation $y^2 = f(x)$ for $E$ minimal at $p$.
\end{remark}

\begin{proposition} \label{PropRamTwoTorsion} Let $p>3$ be a prime and suppose there exists an elliptic curve $E/\QQ$ of conductor $p$ with $E[2](\QQ)=\{\Oc\}$. Then 3 divides one of $h_2(\QQ(\sqrt{p}))$ or $h_2(\QQ(\sqrt{-p}))$.
\end{proposition}
\begin{proof}
Let $\Delta$ be the minimal discriminant of $E$ and let $K=\QQ(\sqrt{\Delta})$. Then, $K$ is one of $\QQ(\sqrt{\pm p})$ or $\QQ(\sqrt{\pm 1})$ and $[L:K]=3$. By Lemma \ref{LemMultRed}, $L/K$ is unramified at $p$, and by standard results it is unramified at primes of good reduction except possibly those above $2$. The existence of such an extension implies that the ray class group of $K$ of modulus $(2)$ has order divisible by $3$ (tame ramification ensures that we can take the modulus $(2)$ rather than $(2^m)$ for $m>1$). This is not the case for $\QQ(\sqrt{\pm 1})$, so $K$ must be one of $\QQ(\sqrt{\pm p})$. 
\end{proof}

\begin{remark}
For $d$ square-free, $d \not \equiv 5 \mod 8$, $3 | h_2(\QQ(\sqrt{d})$ if and only if $3 | h(\QQ(\sqrt{d})$, and for $d<-3$ square-free, $d \equiv 5 \mod 8$, $3$ always divides $h_2(\QQ(\sqrt{d})$ (because 2 is inert and the only units are $\pm 1$). In particular, for primes $p \equiv \pm 1 \mod 8$, Proposition \ref{PropRamTwoTorsion} can be rewritten using only the class group, and for primes $p \equiv 3 \mod 8$, $p>3$, the conclusion of Proposition \ref{PropRamTwoTorsion} always holds.
\end{remark}

We can control the existence of such quadratic fields using the following Cohen-Lenstra type conjecture:

\begin{conj} \label{CLPrecise3Div}
For $a \in \{1,5,7\}$, the density of the set of primes $p \equiv a \mod 8$ such that $3 | h_2(\QQ(\sqrt{p}))$ or $3 | h_2(\QQ(\sqrt{-p}))$ (in the set of all primes $p \equiv a \mod 8$) is:
\begin{itemize}
\item $1 - \prod_{l=1}^\infty (1-3^{-l}) \cdot \prod_{l=2}^\infty (1-3^{-l}) \approx .52939$ for $a \equiv \pm 1 \mod 8$
\item $1 - (\prod_{l=1}^\infty (1-3^{-l}))^2 \approx .68626$ for $a \equiv 5 \mod 8$
\end{itemize}
\end{conj}

\begin{remark}\label{rem:CLWhy}
Conjecture \ref{CLPrecise3Div} follows from the classical Cohen-Lenstra heuristics \cite{CohenLenstra} and the following additional assumptions:
\begin{itemize}
\item The heuristics are valid if, instead of ranging over all fields $\QQ(\sqrt{n})$, we range over fields of the form $\QQ(\sqrt{\epsilon p})$ for $\epsilon=\pm 1$ fixed and $p$ ranging over primes in a fixed congruence class mod $8$.
\item The statistics for imaginary quadratic and real quadratic fields are independent.
\item For real quadratic fields where 2 is inert, the probability that a principal unit generates the unit group mod $(2)$ is $2/3$ (note there are 3 possible elements in $\mathbb{F}_4^\times$, two of which generate it), and this is independent of the probability that 3 divides the order of the class group. 
\end{itemize}

Conjecture \ref{CLPrecise3Div} is compatible with numeric data for $p < 5 \cdot 10^6$. Furthermore, the results of \citep{Bhargava-DH} can be applied to show that a related statistic, the sum of the size of the 3 part of the class group over $\QQ(\sqrt{\epsilon n})$ for a fixed $\epsilon=\pm 1$ and ranging over $n$ square free in a fixed congruence class mod $8$, does behave as if the Cohen-Lenstra heuristics are insensitive to restriction to a conjugacy class. 
\end{remark}

We now give our main result. For $a \in \{ 1,3,5,7 \}$ a congruence class mod $8$, denote by $\P_a$ the set of primes congruent to $a \mod p$. 

\begin{theorem}\label{ThmModDensity} Assuming Conjecture \ref{CLPrecise3Div}, 
\be \limsup_{X\rightarrow{\infty}} \frac{\#\Cond_{P_a}(X)}{\#P_a(X)} \leq
\begin{cases} 1 - \prod_{l=1}^\infty (1-3^{-l}) \cdot \prod_{l=2}^\infty (1-3^{-l}) \approx .52939 & \mbox{if }  a =  1 \textrm{ or } 7 \\
1 - (\prod_{l=1}^\infty (1-3^{-l}))^2 \approx .68626 & \mbox{if } a = 5 
\end{cases}
\ee
\end{theorem}
\begin{proof}
As in Section \ref{subsec-InfPrimeCond} we invoke the result of \label{subsec-InfPC}
Setzer \cite[Thm. 2]{Setzer-ECPC} that for $p\neq2,\,3,\,17$
there is an elliptic curve of conductor $p$ with a rational 2-torsion
point if and only if $p=u^{2}+64$ for some integer value of $u$. Because the set of all integers of the form $u^2+64$ less than $X$ is asymptotically of size $\sqrt{X}$, and by Dirichlet's theorem, for $(a,8)=1$, the number of primes $p\equiv a \mod 8$ less than $X$ is asymptotically $\frac{1}{4}\cdot \log X/X$, Setzer's result implies that the densities of $p\equiv a \mod 8$ such that there exists an elliptic curve of conductor $p$ with non-trivial rational $2$-torsion is $0$. Thus, we need only consider curves with trivial rational $2$-torsion. The result then follows immediately from Proposition \ref{PropRamTwoTorsion} and Conjecture \ref{CLPrecise3Div}.
\end{proof}

Averaging over the different congruence classes, we obtain

\begin{corollary}[Conditional on Conjecture \ref{CLPrecise3Div}]\label{CorTotalDensity}

\be \limsup_{X\rightarrow{\infty}} \frac{\#\Cond_{\P}(X)}{\#P(X)} \leq
1 - (\prod_{l=1}^\infty (1-3^{-l}))^2 \approx .68626. \ee
\end{corollary}

\begin{remark}
If we do not assume independence of real quadratic and imaginary quadratic statistics then we obtain worse but still non-trivial bounds from a weaker version of Conjecture \ref{CLPrecise3Div}. 
\end{remark}

\begin{remark}[Conductors of the form $2p$]\label{Cond2p} If there exists an elliptic curve $E$ of conductor $2p$ with no rational 2-torsion point, then arguing as above we find one of $\QQ(\sqrt{\pm 2p})$ or $\QQ(\sqrt{\pm p})$ has class number divisible by 3 (because the curve has multiplicative reduction at 2 we can apply Lemma \ref{LemMultRed} and thus do not need to worry about the 2-class group even for $p\equiv \pm 3 \mod 8$). If we assume independence of 3-divisibility statistics among these fields and Cohen-Lenstra distributions for fields of the form $\QQ(\sqrt{\pm 2p})$ in addition to the fields $\QQ(\sqrt{\pm p})$, then arguing as in the prime conductor case we obtain (conditional on these assumptions) that the density of primes $p$ such that none of these fields has class number divisible by 3 is $\approx 0.77852$ (which is consistent with numerics for $p < 3 \cdot 10^5$). Thus, conditional on the assumptions, the upper density of $p$ such that there exists an elliptic curve of conductor $2p$ with no rational 2-torsion point is less than $\approx 0.77852$. Using a result of Hadano \cite[Theorem 2]{Hadano-CECRPO2} limiting the 2-torsion on such a curve for $p\equiv \pm 3 \mod 8$ we can remove the assumption on the two-torsion at the price of weakening the bound.
\end{remark}

\section{Conjectures}\label{sec:Conjectures}
\begin{figure}\label{fig:condplot}
\begin{center}
\includegraphics[scale=0.5, natwidth=782, natheight=584]{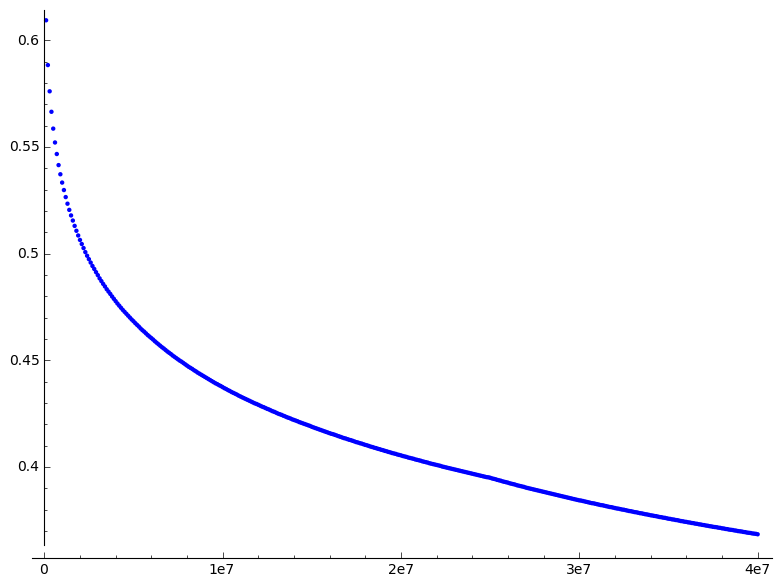}
\caption{Fraction of conductors $\leq X$ for $X\leq 10^8$.}
\end{center}
\end{figure}

In this section state some conjectures which are based on the previous results and on the data presently available. At present, the Cremona \cite{Cremona-ECD} and Stein-Watkins \cite{SteinWatkinsECD} databases provide data on conductors $\leq 10^8$. There are about $61,000$ conductors up to $10^5$ (or about $61\%$) while up to $10^8$ there are only about $31\%$ conductors. The data is plotted in Figure \ref{fig:condplot}.

This leads to the following conjecture:
\bcon\label{con:zero-density}
The set $\Cond$ has density zero in the set of ninth-power free integers, equivalently the set of conductors had density zero in $\N$. Moreover let $\Cond(X)$ be the counting function of $\Cond$. Then 
\be 
\Cond(X)\ll \frac{X}{\log(X)^A},
\ee
for some $A>0$.
\econ

\begin{figure}\label{fig:primeconplot}
\begin{center}
\includegraphics[scale=0.5, natwidth=782, natheight=584]{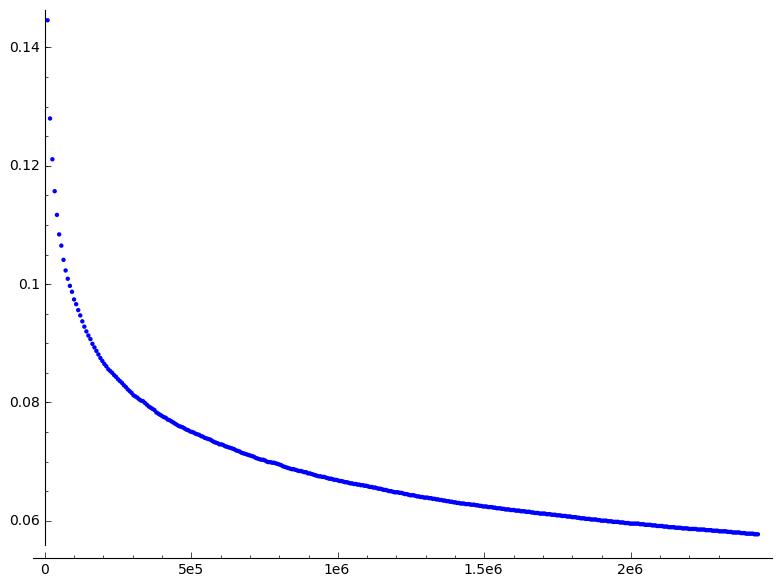}
\caption{Fraction of prime conductors $\leq X$ for $X\leq 10^8$.}
\end{center}
\end{figure}

For prime conductors the numbers here is the corresponding plot (for primes $\leq 10^8$) in \ref{fig:primeconplot}.
This leads us to the following conjecture:
\bcon\label{con:prime-zero-density}
$\Cond_\P$ has density zero in the set of primes.
\econ

Of the $31\%$ conductors $\leq 10^8$, $80\%$ are neither followed by nor are preceded by a conductor. It seems reasonable to expect that almost all conductors are ``loners'' and do not come in clusters of two or more.
\bcon
Almost all conductors $n$ are neither preceded nor followed by a conductor. That is for almost all $n\in\Cond$, neither $n-1$ nor $n+1$ is a conductor.
\econ
The numbers $32,33,34,35,36,37,38,39,$ and $40$ are all conductors. In the Cremona and Stein-Watkins databases one finds there are many consecutive integers which are conductors and also many consecutive integers which are non-conductors. The longest conductor run is of length $18$ and occurs at $n=1130$ so numbers $1130\leq m\leq 1147$ are all conductors of some elliptic curves.  Similarly the data also shows that there are many strings of consecutive non-conductors. The longest consecutive run of non conductors of length  $45$ starts at $n=53649123$. This leads us to the following very optimistic conjecture:

\bcon\label{gap-conjecture}
For any natural number $n\geq2$, there exists at least one chain of consecutive non-conductors of length $\geq n$. In particular non-conductor runs are unbounded.
\econ
This may a be a little too optimistic--given the paucity of the data, but it would very strange to find that consecutive non-conductor runs are bounded.
 
\bibliographystyle{plain}
\bibliography{References}

%
% ----------------------------------------------------------------
%\bibliographystyle{plainnat}
%\bibliography{seanHowe,../../master/master6}
\end{document}